\documentclass[12pt]{article}

\usepackage{geometry,amssymb,amsthm,amsmath,
graphics,enumerate}
\usepackage{times}
\usepackage{amsfonts}
\usepackage{amscd}
\usepackage{url}
  \def\mcolon{\! : \!}

\newcommand\myrestriction{\mathord\restriction}
\def\mr#1{\myrestriction_{#1}}
 \def\mcolon{\! : \!}

\newbox\smilebox
\newbox\anchorbox
\newbox\noanchorbox
\newbox\tempbox

\setbox\smilebox=\hbox{$\smile$}

\def\anchor{\hbox{\vtop{
           \hbox to \wd\smilebox{\hfil\vrule width.4pt height7pt depth1pt\hfil}
           \vskip  -11.5truept
           \hbox to \wd\smilebox{\hfil$\smile$\hfil}}}}
\setbox\anchorbox=\anchor
\def\noanchor{\hbox{\vtop{
           \hbox to \wd\anchorbox{\hfil\anchor\hfil}
           \vskip -14truept
           \hbox to \wd\anchorbox{\hfil/\hfil}}}}
\setbox\noanchorbox=\noanchor

\def\fg#1#2#3{\setbox\tempbox=\hbox{$\scriptstyle{#2}$}
\ifnum\wd\anchorbox>\wd\tempbox\dimen255=\wd\anchorbox
\else\dimen255=\wd\tempbox\fi
{#1\,\vtop{\hbox to \dimen255{\hfil\anchor\hfil}
           \vskip -6truept
           \hbox to \dimen255{\hfil$\scriptstyle{#2}$\hfil}}
           \,#3}}

\def\nfg#1#2#3{\setbox\tempbox=\hbox{$\scriptstyle{#2}$}
\ifnum\wd\noanchorbox>\wd\tempbox\dimen255=\wd\noanchorbox
\else\dimen255=\wd\tempbox\fi
{#1\,\vtop{\hbox to \dimen255{\hfil\noanchor\hfil}
           \vskip -6truept
           \hbox to \dimen255{\hfil$\scriptstyle{#2}$\hfil}}
           \,#3}}

\setbox1=\hbox{$\bot$}

\def\north#1#2{#1\,
\hbox{$\bot$\llap {\hbox to\wd1 {\hfil $/$\hfil}}}
\,#2}

\def\nao#1#2#3{#1\  \hbox{\vtop{
\baselineskip=4pt
\hbox{$\bot$\llap {\hbox to\wd1 {\hfil $/$\hfil}}
\hskip .05em \llap{\hbox{$^{\scriptscriptstyle{a}}$}}}\hbox{$\scriptstyle
{#2}$}}}\, #3}

\def\bp{\par{\bf Proof.}$\ \ $}
\def\abar{\overline{a}}
\def\bbar{\overline{b}}
\def\cbar{\overline{c}}
\def\dbar{\overline{d}}
\def\ebar{\overline{e}}
\def\fbar{\overline{f}}

\def\hbar{\overline{h}}

\def\xbar{\overline{x}}
\def\ybar{\overline{y}}

\def\tp{{\rm tp}}

\def\conc{{\char'136}}
\def\A{{\cal A}}
\def\B{{\cal B}}
\def\C{{\mathfrak  C}}
\def\CC{{\cal C}}
\def\D{{\cal D}}
\def\E{{\cal E}}
\def\S{{\cal  S}}
\def\F{{\cal F}}
\def\FF{{\bf F}}
\def\II{{\bf I}}

\def\Ltau{{\tau_{{\rm ord}}}}

\def\PP{{\mathbb P}}

\def\Q{{\mathbb Q}}
\def\QQ{{\mathbb Q}}

\def\Z{{\mathbb Z}}

\def\tp{{\rm tp}}

\def\Fa0{{\FF^a_{\aleph_0}}}

\def\bp{{\bf Proof.}\quad}

\def\<{\langle}
\def\>{\rangle}

\def\Za0{[Z]^{\aleph_0}}

\def\U{{\cal U}}
\def\code#1{{\ulcorner #1 \urcorner}}
\def\At{{\bf At_T}}

\def\pcl{{\rm pcl}}
\newtheorem{Theorem}{Theorem}[subsection]
\newtheorem{Claim}[Theorem]{Claim}

\newtheorem{Definition}[Theorem]{Definition}
\newtheorem{Notation}[Theorem]{Notation}

\newtheorem{Example}[Theorem]{Example}

\newtheorem{Lemma}[Theorem]{Lemma}

\newtheorem{Fact}[Theorem]{Fact}

\newtheorem{Construction}[Theorem]{Construction}

\def\Gdot{{\bf \dot{G}}}

\begin{document}


\title{Constructing many atomic models in $\aleph_1$}

\author{John T. Baldwin\thanks{Research partially supported by Simons
travel grant G5402}
\\University of Illinois at Chicago\\
\and Michael C. Laskowski\thanks{Partially supported by NSF grant
DMS-1308546}\\University of Maryland\\
\and Saharon Shelah\thanks{This research was partially supported by
NSF grant DMS 1101597. The third author was partially supported by
the
European Research Council grant 338821. This is paper 1037 in the Shelah Archive}\\
 Hebrew University of Jerusalem and Rutgers University }

\maketitle

\section{Introduction}


 As has been known since at least \cite{Sh87a} and is
carefully spelled out in Chapter 6 of \cite{Baldwincatmon}, for every
complete sentence $\psi$ of $L_{\omega_1,\omega}$ (in a countable
vocabulary $\tau$) there is a complete, first order theory $T$ (in a
countable vocabulary extending $\tau$) such that the models of $\psi$
are exactly the $\tau$-reducts of  the atomic models of $T$.
  This
paper is written entirely in terms of the class $\At$ of atomic
models of a complete first order theory $T$, but applies to
$L_{\omega_1,\omega}$ by this translation.

Our main theorem, Theorem~\ref{big}, asserts: Let $T$ be any complete
first-order theory in a countable language with an atomic model. If
the pseudo-minimal types are not dense, then there are $2^{\aleph_1}$
pairwise non-isomorphic, full\footnote{An atomic model $M$ is {\em
full} if $|\phi(M,\abar)|=||M||$ for every non-pseudo-algebraic
formula $\phi(x,\abar)$  (See Definition~\ref{arbitrary}.) with
$\abar$ from $M$.} atomic models of $T$, each of size $\aleph_1$.

 The first section states some old observations about atomic models and
develops a notion of `algebraicity', dubbed pseudo-algebraicity for
clarity, that is relevant in this context.  We introduce the relevant
analogue to strong minimality, pseudo-minimality, and state the
pseudo-minimals dense/many models dichotomy. Section~\ref{2} expounds
a transfer technique, already used in \cite{BL} and \cite{BLS} and
 applied here prove to Theorem~\ref{big}.  The gist of the method
is to prove a model theoretic property is consistent  with ZFC by
forcing and then extend the model $M$ of set theory witnessing this
result to a model $N$, preserving the property and such that  the
property is absolute between $V$ and $N$. Section~\ref{3} describes a
forcing construction, which together with the results of
Section~\ref{2}, yields a proof of Theorem~\ref{big} in
Section~\ref{5}.

The authors are grateful to Paul Larson and Martin Koerwien for many insightful conversations.

\section{A notion of algebraicity}
\label{1} \numberwithin{Theorem}{section} \setcounter{Theorem}{0}

Throughout this paper, $T$ will always denote a complete, first-order theory in a countable language that has an atomic model.   By definition, a model $M$ of $T$ is atomic
 if every finite tuple $\abar$ from $M$ realizes a complete formula\footnote{Recall
 that
   $\phi(\xbar)$ is a {\em complete formula in T} if $\phi(\xbar)$ is
the generator of a principal type, i.e. for every $\psi(\xbar)$,
$T\vdash (\forall \xbar)[\phi(\xbar) \rightarrow \psi(\xbar)]$ or
$T\vdash (\forall \xbar)[\phi(\xbar) \rightarrow \neg\psi(\xbar)]$ .
}. The existence of an atomic model is equivalent to the statement
that `every consistent formula $\phi(\xbar)$ has a complete formula
$\psi(\xbar)$ implies it.' Equivalently, $T$ has an atomic model if
and only if, for every $n\ge 1$, the isolated complete $n$-types are
dense in the Stone space $S_n(\emptyset)$. We recall some old results
of Vaught concerning this context.

\begin{Fact} \label{Vaught} Let $T$ be any complete theory in a countable language having an atomic model.
Then:
\begin{enumerate}
\item  $\At$ is $\aleph_0$-categorical, i.e., every pair of countable atomic models are isomorphic;
\item  $\At$ contains an uncountable model if and only if some/every countable model of $\At$ has a proper elementary extension.
\end{enumerate}
\end{Fact}

 The only known arguments for proving amalgamation and thus constructing monster models for
 $\At$ invoke the continuum hypothesis and so are not useful for our
 purposes. Nevertheless, we argue that many concepts of interest are in fact model
 independent.



In first-order model theory, if a formula $\phi(x,\abar)$ is
algebraic, then its solution set cannot be increased in any
elementary extension, i.e., if $\abar\subseteq M\preceq N$, then
$\phi(M,\abar)=\phi(N,\abar)$.  However, in the atomic case, the
analogous phenomenon can be witnessed by non-algebraic formulas.  For
example, $(\Z,S)$, the integers with a successor function, is an
atomic model of its theory.  The formula `$x=x$' is not algebraic,
yet $(\Z,S)$ has no proper {\bf atomic} elementary extensions. This
inspires the following definition:

%


\begin{Definition}  \label{countable}  {\em
Let $M\in\At$ be countable\footnote{In
    Definition~\ref{countable} it would be equivalent to
    restrict to countable and $M$ and allow arbitrary cardinality
    for $N$.  It would {\bf not} be equivalent to assert for
    arbitrary $M$: ``$\phi(x,\abar)$ is pseudo-algebraic in $M$
    if and only if $\phi(M,\abar)=\phi(N,\abar)$ for every
    $N\succeq M$."  To see the distinction, consider the extreme
    case where $M$ is an uncountable atomic model that is
    maximal, i.e., has no proper atomic elementary extension.}. A formula $\phi(x,\abar)$ is {\em
pseudo-algebraic in $M$} if $\abar$ is from $M$, and
$\phi(N,\abar)=\phi(M,\abar)$ for every countable $N\in\At$ with
$N\succeq M$. }
\end{Definition}



 The strong $\aleph_0$-homogeneity (any two finite sequences realizing the same type over the emptyset are automorphic)
  of the countable atomic
model of $T$ yields immediately that pseudo-algebraicity  truly
depends only on the type of $\abar$ over the emptyset.  That is, if
$M,M'\in\At$ are each countable and $\tp(\abar,M)=\tp(\abar',M')$,
then $\phi(x,\abar)$ is pseudo-algebraic in $M$ if and only if
$\phi(x,\abar')$ is pseudo-algebraic in $M'$. This observation allows
us to extend the notion of pseudo-algebraicity to arbitrary atomic
models of $T$.

\begin{Definition}  \label{arbitrary}  {\em  Let $N\in\At$ have arbitrary cardinality.
\begin{enumerate}

\item A formula $\phi(x,\abar)$ is {\em pseudo-algebraic in $N$}
    if $\abar$ is from $N$, and $\phi(x,\abar)$ is
    pseudo-algebraic in $M$ for some (equivalently, for every)
    countable $M\preceq N$ containing $\abar$.

\item  An element $b\in N$ is {\em pseudo-algebraic over
    $\abar$ inside $N$}, written $b\in\pcl(\abar,N)$, if $\tp(b/\abar,N)$
contains a formula that is pseudo-algebraic in $N$.
\item  Given an infinite subset $A\subseteq N$, {\em $b$ is
    pseudo-algebraic over $A$ in $N$}, written   $b\in\pcl(A,N)$,
    if and only if $b\in\pcl(\abar,N)$ for some finite $\abar\in
    A^n$.
\end{enumerate}
}
\end{Definition}



As the language of $T$ is countable,
for any complete formula $\theta(\ybar)$, there is a formula
$\psi(x,\ybar)$ of $L_{\omega_1,\omega}$ such that $T\cup\{\psi(x,\ybar)\}\vdash\theta(\ybar)$
and for every atomic
$M$, every $\abar\in\theta(M)$, and every $b\in M$:
$$b\in\pcl(\abar,M)\qquad\hbox{if and only if}\qquad M\models\psi(b,\abar)$$

Note that this notion allows us to reword Fact~\ref{Vaught}(2):  $T$
has an uncountable atomic model if and only if `$x=x$' is not
pseudo-algbraic. Here is a second example.

\begin{Example}
{\em
Let $L=\{A,B,\pi,S\}$ and $T$ say that $A$ and $B$ partition
	the universe with $B$ infinite,
	$\pi:A\rightarrow B$ is a total surjective
	function and $S$ is a successor function on $A$ such that
	every $\pi$-fiber is the union of $S$-components. 	
	A model $M\models T$  is atomic  if every
	$\pi$-fiber contains exactly one $S$-component.
	Now choose elements $a,b\in M$ for such an $M$ such that
	$a\in A$ and $b\in B$ and $\pi(a) =b$. Clearly, $a$ is not
algebraic	over $b$ in the classical sense, but $a\in \pcl(b,M)$.
}
\end{Example}

Recall that a {\em $t$-construction over $B$} is a sequence $\langle a_i:i<
\omega\rangle$ such that, letting $A_i$ denote $B\cup\{a_j:j<i\}$,
$\tp(a_i/A_i)$ is generated by a complete formula.

The notion of pseudo-algebraicity has many equivalents.  Here are
some we use below.
\begin{Lemma}\label{pcleq}  Suppose $M\in\At$ and $b,\abar$ are from $M$.  The following are equivalent:
\begin{enumerate}
\item  $b\in\pcl(\abar,M)$;
\item  For every   $N\preceq M$, if $\abar\in N^n$, then $b\in N$;
\item  $b$ is contained inside any maximal $t$-construction sequence $\<a_\alpha:\alpha<\beta\>$ over $\abar$ inside $M$.
\end{enumerate}
\end{Lemma}

For (3) note that as $T$ has an atomic model, a maximal
$t$-construction sequence over a finite set is the universe of a
model.
%

Here is one application of Lemma~\ref{pcleq}.

\begin{Lemma}  \label{outside}
Suppose that $M\in\At$, $\abar$ is from $M$, but $\phi(x,\abar)$ is not pseudo-algebraic in $M$.  Then for every finite $\ebar$ from $M$,
there is $b\in \phi(M,\abar)$ with $b\not\in \pcl(\ebar,M)$.
\end{Lemma}

\bp  We may assume $\abar\subseteq \ebar$. Choose a countable
$M^*\preceq M$ containing $\ebar$ and, by non-pseudo-algebraicity and
Definition~\ref{countable}, choose  a countable $N^*\in\At$  with $N^*\succeq M^*$ and $b^*\in
\phi(N^*,\abar)\setminus \phi(M^*,\abar)$.
%
As $N^*$ is
countable and atomic, choose an elementary embedding
$f:N^*\rightarrow M$ that fixes $\ebar$ pointwise.  Then
$f(b^*)\in\phi(M,\abar)$
 and $f(b^*)\not\in\pcl(\ebar,M)$ as witnessed by $f(M^*)$ and Lemma~\ref{pcleq}(2).
$\qed_{\ref{outside}}$


\medskip

In general, the notion of pseudo-algebraic closure gives rise to a
reasonable closure relation.  All of the standard van der Waerden
axioms for a dependence relation hold in general, with the exception
of the Exchange Axiom.  Our next definition isolates  those formulas
on which exchange (and a bit more) hold.

\begin{Definition}\label{psmdef} {\em   Let $M$ be any atomic model and let $\abar$ be from $M$.
\begin{itemize}
\item  A complete formula $\phi(x,\abar)$ is {\em pseudo-minimal}
    if it is not pseudo-algebraic, but for every
    $\abar^*\supseteq\abar$ and $c$ from $M$ and for every
    $b\in\phi(M,\abar)$, if $c\in\pcl(\abar^*b,M)$ but
    $c\not\in\pcl(\abar^*,M)$, then $b\in\pcl(\abar^*c,M)$.
\item  The class $\At$ has {\em density of pseudo-minimal types}
    if for some/every $M\in\At$,  for every
     non-pseudo-algebraic formula $\phi(x,\abar)$, there is
    $\abar^*\supseteq\abar$ from $M$ and a pseudo-minimal formula
 $\psi(x,\abar^*)$ such that
 $\psi(x,\abar^*)\vdash\phi(x,\abar)$.
\end{itemize}
}
\end{Definition}


It is immediate that if there is a non-pseudo-algebraic formula then
$T$ has an atomic model in $\aleph_1$, so also if pseudo-minimal types
are not dense,  then $T$ has an atomic model in $\aleph_1$.  The main
Theorem of this paper is the following:

\begin{Theorem}  \label{big}
Let $T$ be any complete first-order theory in a countable language with an atomic model.
If the pseudo-minimal types are not dense, then there are $2^{\aleph_1}$ pairwise non-isomorphic,
full, atomic models
of $T$, each of size $\aleph_1$.
\end{Theorem}

\section{A technique for producing many models of power $\aleph_1$}  \label{2}
\numberwithin{Theorem}{subsection} \setcounter{Theorem}{0}

 The
objective of this section is to prove the transfer
Theorem~\ref{transfer1} that allows the construction  (in ZFC) of
many atomic models of a first order theory $T$ in two steps. First
force to find a model $(M,E)$ of set theory in which a model of $T$
is coded by stationary sets. Then apply the transfer theorem to code
a family of such models in ZFC.

The method expounded here has many precursors. Among the earliest are
the treatment of Skolem ultrapowers in \cite{Keislerbook} and the
study of elementary extensions of  models of set theory in \cite{KeislerMorley} and
\cite{Hutch}. Paul Larson introduced the use of iterated generic
ultrapowers (used in the different context of Woodin's $\PP$-max
forcing) in a large cardinal context in \cite{FL, FKLM} and the
general method is abstracted in \cite{Larsonit}. The model theoretic
technique used here is described in \cite{BL} and \cite{BLS}. We
formulate a general metatheorem for the construction.

 The first subsection
describes how to define and maintain satisfaction of formulas in a pre-determined,
countable fragment $L_\A$ under elementary extensions of
$\omega$-models of set theory.
 Most of this is well-known; we emphasize that only an $\omega$-model
 and not transitivity is necessary to correctly code sentences of
 $L_{\omega_1,\omega}$.
The second subsection surveys known results about
$M$-normal ultrapowers, and Theorem~\ref{transfer1} is proved in the
third subsection.

\subsection{Coding $\tau$-structures into non-transitive models of set theory}

In this section, we fix an explicit encoding of a pre-determined
countable fragment $L_\A = L_\A(\tau)$ of $L_{\omega_1,\omega}(\tau)$
for a     countable vocabulary $\tau$ into an $\omega$-model $(M,E)$
satisfying $ZFC$. The specific form of this encoding is not
important, but it is useful for the reader to see what we assume
about $M$ in order that satisfaction is computed `correctly' for
every formula of $L_\A$. It will turn out that everything works
wonderfully (even when $(M,E)$ is non-transitive) provided $(M,E)$ is
an $\omega$-model (that is $\omega^M=  \omega^V$), because this
guarantees a formula of $L_\A$ does not gain additional conjuncts or
disjuncts in an elementary extension that is also  an $\omega$-model.



\begin{Definition} {\em  We say $(M,E)$ is an {\em $\omega$-model of set theory} if $(M,E)\models ZFC$,
$(\omega+1)^{M,E}=\omega+1$, and for $n,m\in\omega+1$, $(M,E)\models nEm$ if and only if $n\in m$.
}
\end{Definition}

Fix any countable vocabulary (sometimes called language) $\tau$.  In
what follows, we will assume that $\tau$ is relational with
$\aleph_0$ $n$-ary relation symbols $R^n_m$, but the generalization
to other countable languages is obvious.

%

\begin{Definition} {\em  Fix a particular countable fragment $L_\A = L_\A(\tau)$ of $L_{\omega_1,\omega}(\tau)$.
\begin{itemize}
\item  A {\em Basic G\"odel number} has the form $\<0,n,m\>$,
    where $n,m\in\omega$. \ We write this as $\code{R^n_m}$.
\item Let $BG_\tau$ denote the set of Basic G\"odel numbers.  We
    now define by induction the set $G_{L_\A}$  of G\"odel
    numbers of $L_\A$-formulas.
    \begin{enumerate}
\item  $\code{v_i}=\<1,i\>$;
\item $\code{R^n_m(v_{i_1}, \ldots v_{i_n})}=\<\code{R^n_m},
    \code{v_{i_1}},\dots ,\code{v_{i_n}}\>$
\item  $\code{\phi=\psi}=\<2,\code{\phi},\code{\psi}\>$;
\item  $\code{\phi\wedge\psi}=\<3,\code{\phi},\code{\psi}\>$;
\item  $\code{\exists v_i\phi}=\<4,\code{v_i},\code{\phi}\>$;
\item  $\code{\neg\phi}=\<5,\code{\phi}\>$;
\item  If $\psi=\bigwedge_{i\in\omega}\theta_i$ and $\psi\in
    L_\A$, then $\code{\psi}=\<6,f_\psi\>$, where $f_\psi$ is the
    function with domain $\omega$ and
    $f_\psi(i)=\code{\theta_i}$.
    \end{enumerate}

\end{itemize}
}
\end{Definition}

\begin{Definition}  {\em  For a given countable fragment $L_\A$, we say an $\omega$-model $(M,E)$ {\em supports $L_\A$}
if  $G_{L_\A}\in M$ and $G_{L_\A}\subseteq M$.
}
\end{Definition}


Note that $BG_{\tau}$ and $G_{L_\A}$  are defined in $V$ but they are
correctly identified by an $(M,E)$ that supports $L_\A$. More
precisely, the following lemma is immediate.

%

\begin{Lemma}  If $(M,E)$ is an $\omega$-model of set theory supporting $L_\A$, then
both $BG_\tau$ and $G_{L_\A}$ are definable subsets of $M$. Furthermore,
if $(N,E)\succeq (M,E)$ is also an $\omega$-model, then
$BG_\tau^{N,E}=BG_\tau^{M,E}$, $(N,E)$ supports $L_\A$,
$G_{L_\A}^{N,E}=G_{L_\A}^{M,E}$, and
$\code{\phi}^{N,E}=\code{\phi}^{M,E}$ for every $\phi\in L_\A$.
\end{Lemma}

\begin{Definition}
{\em  Suppose $(M,E)$ is an $\omega$-model of set theory, and we have
fixed a countable vocabulary $\tau$. A $\tau$-structure $\B=(B,\dots)$
is {\em inside $(M,E)$ via $g$} if the universe $B\in M$, $g\in M$ is
a function with domain $BG_\tau\cup\{\emptyset\}$, $g(\emptyset)=B$ and
for each $(n,m)\in\omega^2$, $g(\code{R^n_m})=R^n_m(\B)$. }
\end{Definition}

\begin{Definition} \label{extension}
{\em If $(M,E)$ is an $\omega$-model of set theory, a
$\tau$-structure $\B$ is inside $(M,E)$ via $g$, and
$(N,E)\succeq(M,E)$ is an $\omega$-model, then $\B^N$ denotes the
$\def\PP{{\mathbb P}}$-structure with universe $g(\emptyset)^N$ and
relations $R^n_m(\B^N)=g(\code{R^n_m})^N$. }
\end{Definition}

Clearly, $\B^N$ is inside $(N,E)$ via $g^N$. Again using the fact
that we are working with $\omega$-models, the following is immediate.

\begin{Lemma}\label{embed}  Suppose $(M,E)$ is an $\omega$-model of set theory supporting $L_\A$
 and a $\tau$-structure $\B$ is inside  $(M,E)$ via $g$.
Then there is a unique $h\in M$, $h\mcolon G_{L_\A}\rightarrow M$
extending $g$ such that $h(\code{\psi})=\psi(\B)$ for every $\psi\in
L_\A$.
\end{Lemma}

\subsection{$M$-normal ultrapowers}

The idea of using $M$-normal ultrafilters to construct many elementary chains of models of set theory is not new, and the definitions and results of this subsection
are presented here for the
convenience of the reader.

Fix a countable $\omega$-model $(M,E)$ of set theory.
Since $M$ is countable, so is the set $\omega_1^M$.
As notation, let $$\C=\{B\subseteq\omega_1^M:M\models `B\ \hbox{is club'}\}$$
In what follows, a function $f$ with domain $\omega_1^M$ is {\em regressive} if $f(\alpha)<\alpha$ for all $\alpha>0$.

\begin{Definition} {\em
An {\em $M$-normal ultrafilter}
 $\U$ is an ultrafilter on the set $\omega_1^M$ such that
\begin{itemize}
\item  $\C\subseteq\U$; and
\item  For every regressive $f:\omega_1^M\rightarrow\omega_1^M$ with
$f\in M$, $f^{-1}(\beta)\in\U$ for some $\beta\in\omega_1^M$.
\end{itemize}
}
\end{Definition}

We record an Existence Lemma for $M$-normal ultrafilters.

\begin{Lemma} \label{existence} Suppose $A\subseteq\omega_1^M$ and $A\in M$.
Then there is an $M$-normal ultrafilter $\U$ with $A\in\U$ if and only if
$M\models `A\ \hbox{is stationary'}$.
\end{Lemma}

\bp  Clearly, if
$M\models `A\ \hbox{is non-stationary'}$, then there is some $B\in\C$ such that $A\cap B=\emptyset$, so no $M$-normal
ultrafilter can contain $A$.  For the converse, enumerate the regressive functions in $M$ by $\<f_n:n\in\omega\>$.
We construct a nested, decreasing sequence $\<A_n:n\in\omega\>$ of subsets of $\omega_1^M$ such that each $A_n\in M$ and
$M\models `A_n\ \hbox{is stationary'}$ as follows:
Put $A_0:=A$ and given $A_n$, by Fodor's Lemma (in $M$!) choose a stationary $A_{n+1}\subseteq A_n$ and $\beta_n$ such that
$f_n[A_{n+1}]=\{\beta_n\}$.

As $\C\cup\{A_n:n\in\omega\}$ has f.i.p., (now working in $V$) it follows that there is an ultrafilter $\U$ containing these sets.
Any such $\U$ must be $M$-normal.
$\qed_{\ref{existence}}$

We record three consequences of $M$-normality.

\begin{Lemma}  \label{normal}
Suppose that $\U$ is an $M$-normal ultrafilter on $\omega_1^M$.  Then:
\begin{enumerate}
\item  If $A\in\U\cap M$, then $M\models `A\ \hbox{is stationary'}$;
\item  If $A\in\U\cap M$, $f\in M$, and $f:A\rightarrow \omega_1^M$ is regressive, then $f^{-1}(\beta)\in\U$ for some $\beta\in\omega_1^M$; and
\item  If
$\<A_n:n\in\omega\>\in M$ and every $A_n\in \U\cap M$, then $A=\bigcap_{n\in\omega} A_n\in \U\cap M$.
\end{enumerate}
\end{Lemma}

\bp (1)  Choose $A\in\U\cap M$.  To see that $A$ is stationary in $M$,
choose any $B\in M$ such that $M\models `B\ \hbox{is club'}$.  Then $B\in\C\subseteq\U$.  As $\U$ is a proper filter, $A\cap B$ is non-empty.

(2)  This is `completely obvious' but rather cumbersome to prove
precisely.


Given $f:A\rightarrow\omega_1^M$, by intersecting with the club
$D:=\omega_1^M\setminus\omega$, we may assume $A\subseteq D$.
Define $g:\omega_1^M\rightarrow\omega^M_1$ by

$$g(\delta)=
 \left\{\begin{array}{ll}
f(\delta) &\mbox{if $\delta\in A$ and $f(\delta)\ge\omega$}\\
f(\delta)+1 &\mbox{if $\delta\in A$ and $f(\delta)<\omega$}\\
0 & \mbox{if $\delta\not\in A$}
\end{array}
\right. $$
Then $g\in M$ and $g$ is regressive, hence $g^{-1}(\beta)\in\U$ for some $\beta$.
As $g^{-1}(0)$ is disjoint from $A$ and $A\in\U$, $\beta\neq 0$.  Thus, $g^{-1}(\beta)\subseteq A$.
It follows that either $f^{-1}(\beta)\in \U$ (when $\beta\ge\omega$) or $f^{-1}(\beta-1)\in\U$ (when $\beta<\omega$).

(3)  Assume not.  Let $B:=\omega_1^M\setminus A\in \U\cap M$.  As in (2) we may assume $B\subseteq (\omega_1^M\setminus\omega)$.
Define $f:B\rightarrow\omega$ by
$$f(\delta)=\hbox{least $n$ such that $\delta\not\in A_n$}$$
As $f$ is regressive, we get a contradiction from (2).
$\qed_{\ref{normal}}$

%

\medskip\par

Given $M$ and an $M$-normal ultrafilter $\U$, we form
the ultraproduct $Ult(M,\U)$ as follows:

First, consider the (countable!) set of functions $f:\omega_1^M\rightarrow M$
with $f\in M$.  There is a natural equivalence relation $\sim_\U$ defined by
$$f\sim_\U g \quad\Leftrightarrow\quad \{\delta\in\omega_1^M:f(\delta)=g(\delta)\}\in\U$$
The objects of $Ult(M,\U)$ are the equivalence classes $[f]_\U$, and
we put
$$Ult(M,\U) \models [f]_\U E [g]_\U\quad\Leftrightarrow\quad \{\delta\in\omega_1^M:
f(\delta) E g(\delta)\}\in\U.$$
%

%

For each $a\in M$, we have the constant function
$f_a:\omega_1^M\rightarrow M$ defined by $f_a(\delta)=a$ for every
$\delta\in\omega_1^M$.  Every such function $f_a\in M$, hence we get
an embedding
$$j:M\rightarrow Ult(M,\U)$$
defined by $j(a)=[f_a]_\U$.

The following Lemmas summarize the results we need.

\begin{Lemma}\label{basicnormal}  Suppose that $(M,E)$ is  a countable $\omega$-model of set theory and
$\U$
is any $M$-normal ultrafilter on $\omega_1^M$.
Then:
\begin{enumerate}
\item  $N:=Ult(M,\U)$ is a countable $\omega$-model and
    $j:(M,E)\rightarrow (N,E)$ is elementary.



\item If $a\in M$ and $M\models `a\ \hbox{is countable'}$ then
    $j(a) =j[a]=_{df} \{j(x):x E a\}$.

\item The image $j[\omega_1^M]=_{df}\{j(a):a\in\omega_1^M\}$ is a
    proper initial segment of $\omega_1^N$ with $[id]_\U$ the
    least element of $\omega_1^N\setminus j[\omega_1^M]$.

%

\end{enumerate}
\end{Lemma}

\bp We begin with (2).  Fix $a\in M$ with $M\models `a\ \hbox{is
countable'}$ and abbreviate $M \models aEb$ by $aEb$.  First, for
every $b E a$, $f_b(\delta) E f_a(\delta)$ for every $\delta E
\omega_1^M$, so $j(b) E j(a)$ by \L o\'s's theorem. Conversely, to
show $j(a) \subseteq j[a]$, choose any $g:\omega_1^M\rightarrow M$
with $g\in M$ such that $[g]_\U\neq [f_b]_\U$ for every $b E a$.
Towards showing that $[g]_\U\neg E j(a)$, choose, using the
countability of $a$ in $M$, a surjection $\Phi:\omega\rightarrow a$
with $\Phi\in M$. In $M$, let
$$A_n=\{\delta E \omega_1^M:g(\delta)\neq \Phi(n)\}.$$ By
separation, each $A_n \in M$ and recursion, since $M$ is an
$\omega$-model,
 $\<A_n:n\in\omega\>\in M$ and each $A_n\in\U\cap M$. Thus, by
Lemma~\ref{normal}(3), $A:=\bigcap_{n\in\omega} A_n\in\U\cap M$.
Since $g(\delta)\neg E a$ for every $\delta\in A$, the fact that
$A\in \U$ implies that $[g]_\U\neg E j(a)$.


%

As for (1), that $j:(M,E)\rightarrow (N,E)$ is elementary is  the \L
o\'s theorem. $N$ is clearly countable, as there are only countably
many functions in $M$, and it is an $\omega$-model by (2). As for
(3), that $j(\omega_1^M)$ is an initial segment of $\omega_1^N$
follows from (2), and the minimality of $[id]_\U$ in the difference
follows from Fodor's Lemma in $M$. $\qed_{\ref{basicnormal}}$


We now drop the pedantry of keeping exact track of the embedding $j$
and just write $M \preceq N$.

\begin{Lemma}\label{Lnormal}  Suppose that $(M,E)$ is  a countable $\omega$-model of set theory that supports $L_\A$
and
let $\B=(B,\dots)$ be an $L$-structure inside $(M,E)$ via $g$.  Given any
 $M$-normal ultrafilter $\U$ on $\omega_1^M$, let
$N=Ult(M,\U)$ and let $\B^N$ be the $L$-structure formed as in
Definition~\ref{extension} with $h$ as in Lemma~\ref{embed}. Then:
\begin{enumerate}
\item  For every $L_\A$-formula $\psi(x_1,\dots,x_n)$ and all $[f_1]_\U,\dots,[f_n]_\U$ with each $f_i:\omega_1^M\rightarrow B$,
$$\B^N\models\psi([f_1],\dots,[f_n])\Longleftrightarrow
\{\alpha\in\omega_1^M:(f_1(\alpha),\dots,f_n(\alpha))\in h(\code{\psi})\}\in \U$$
\item  The induced embedding $j:\B\rightarrow\B^N$ is $L_\A$-elementary; and
\item  If $\omega_1^M\subseteq B$ and $\theta(x)\in L_\A$ has one free variable, then
$\B^N\models\theta([id]_\U)$ if and only if $\{\alpha\in\omega_1^M:\alpha\in h(\code{\theta})\}\in\U$.
\end{enumerate}
\end{Lemma}

\subsection{A transfer theorem}
We bring together the methods of the previous subsections into a
general transfer theorem.  Recall that we are using Roman letters (M)
for models of set theory, Gothic ($\B$) for $\tau$-structures and
$\B^M$ denotes a structure supported in $M$, and for a
$\tau$-relation $P$, $P^{\B}$ denotes the elements of $\B$ satisfying
$P$.

\begin{Theorem}\label{transfer1}
Fix a vocabulary $\tau$ with a distinguished unary predicate $P$ and
fix a countable fragment $L_\A= L_\A(\tau)\subset
L_{\omega_1,\omega}(\tau)$. SUPPOSE there is a countable,
$\omega$-model $(M,E)$ of set theory supporting $L_\A$ and there is a
$\tau$-structure $\B=(B,\dots)$ inside $M$ via $g$ satisfying:
\begin{itemize}
\item  $P^\B\subseteq\omega_1^M\subseteq B$;
\item $M\models `P^\B\ \hbox{is stationary/costationary'}$.
\end{itemize}

THEN for every $X\subseteq\omega_1$ (in $V$!) there is an $\omega$-model $(N_X,E)\succeq (M,E)$
and a continuous, strictly increasing\footnote{The function $t_X$ need not be an element of $N_X$.}
$t_X:\omega_1\rightarrow\omega_1^{N_X}$
satisfying:
\begin{itemize}
\item  $|N_X|=\aleph_1$ and $(\omega_1^{N_X},E)$ is an $\aleph_1$-like linear order;
\item  for all $\alpha\in\omega_1$, $\B^{N_X}\models P(t_X(\alpha))$ if and only if $\alpha\in X$.
\end{itemize}
\end{Theorem}

\bp Fix any $X\subseteq \omega_1$. We construct a continuous chain
$\<M_\alpha:\alpha\in\omega_1\>$ of $\omega$-models of set theory as
follows: Put $M_0:=(M,E)$ and at countable limit ordinals, take
unions.  Now suppose $M_\alpha$ is given.  Choose an
$M_\alpha$-normal ultrafilter $\U_\alpha$ such that
$P^{M_\alpha}\in\U_\alpha$ if and only if $\alpha\in X$.  The
existence of such a $\U$ follows from Lemma~\ref{existence}, since by
elementarity, letting $\B_\alpha$ denote $\B^{M_\alpha}$, we have
that
$$M_\alpha\models `P^{\B_\alpha}\ \hbox{is a stationary/costationary subset of $\omega_1$'}$$

Given such a chain, put $N_X:=\bigcup\{M_\alpha:\alpha\in \omega_1\}$ and define $t_X:\omega_1\rightarrow\omega_1^{N_X}$
by $t_X(\alpha)=[id]_{\U_\alpha}$.
$\qed_{\ref{transfer1}}$

This result extends easily to $L(Q)$ and the somewhat more
complicated version for $L(aa)$ is treated in section 2 of \cite{BL}.

\section{The relevant forcing}  \label{3}

Throughout this section, we have a fixed atomic class $\At$ that
contains uncountable models, for which the pseudo-minimal types are
not dense. The objective of this section is introduce a class of
$I^*$ of expansions of linear orders, develop the notion of a model
$N\in\At$ being {\em striated} by such an order, and prove
Theorem~\ref{forcing}, which uses the failure of density of
pseudo-minimal types to force the existence of a striated model
capable of encoding a nearly arbitrary subset of $\omega_1$.

\subsection{A class of linear orders}

Recall that a linear order is $\aleph_1$-like if every initial
segment is countable.  It is well-known that there are $2^{\aleph_1}$
$\aleph_1$-like linear orders of cardinality $\aleph_1$.  An
accessible account of this proof, which underlies this entire paper,
appears on page 203 of \cite{Markerbook}.  The key idea of that
argument is to code a stationary set of cuts which have a least upper
bound. In the current paper, the coding is not so sharp. Instead, we
force an atomic model of $T$ that codes a stationary set by
infinitary formulas defined using $\pcl$.

We begin by describing a class of $\aleph_1$-like linear orders, colored by a unary predicate $P$ and an
equivalence relation $E$ with convex classes.
This subsection makes no reference to the class $\At$.

\begin{Definition} {\em  Let $\Ltau=\{<,P,E\}$ and let $\II^*$ denote the collection of $\Ltau$-structures $(I,<,P,E)$ satisfying:
\begin{enumerate}
\item $(I,<)$ is an $\aleph_1$-like dense linear order with minimum element $\min(I)$ (i.e., $|I|=\aleph_1$, but $pred_I(a)$ is countable for every $a\in I$);
\item $P$ is a unary predicate and $\neg P(\min(I))$;
\item  $E$ is an equivalence relation on $I$ with convex classes such that
\begin{enumerate}
\item  If $t=\min(I)$ or if $P(t)$ holds, then $t/E=\{t\}$;
\item  Otherwise, $t/E$ is a (countable) dense linear order
    without endpoints.
\end{enumerate}
\item  The quotient $I/E$ is a dense linear order with minimum
    element, no maximum element, such that both sets
    $\{t/E:P(t)\}$ and $\{t/E:\neg P(t)\}$ are dense in it.
\end{enumerate}
}
\end{Definition}

Note that for $s\in I$, we denote the equivalence class of $s$ by
$s/E$ and the predecessors of the class by $<s/E$. We are interested
in well-behaved proper initial segments $J$ of orders $I$ in $\II^*$.

\begin{Definition} \label{suitdef} {\em  Fix $(I,<,P,E)\in\II^*$.  A proper initial segment $J\subseteq I$ is {\em suitable} if, for every $s\in J$ there is $t\in J$, $t>s$,
with $\neg E(s,t)$.
}
\end{Definition}

Note that if $J\subseteq I$ is suitable, then $J$ is a union of $E$-classes and that there is no largest $E$-class in $J$.  Accordingly, there are three possibilities for
$I\setminus J$:
\begin{itemize}
\item  $I\setminus J$ has a minimum element $t$.  In this case, it must be that
$t/E=\{t\}$.
\item  $I\setminus J$ has no minimum $E$-class.  In this case, we call $J$ {\em seamless}.
\item  $I\setminus J$ has a minimum $E$-class that is infinite.  This will be our least interesting case.
\end{itemize}

We record one easy Lemma.


\begin{Lemma} \label{seamless}
 If $(I,<,P,E)\in\II^*$ and $J\subseteq I$ is a seamless proper initial segment, then for every finite
  $\S\subseteq I$ and $w\in J$ such that    
$w > \S\cap J$, there is an automorphism $\pi$ of $(I,<,P,E)$ that
fixes $\S$ pointwise, and $\pi(w)\not\in J$.
\end{Lemma}

\bp  Fix $I,J,\S$ as above.  As $J$ is seamless, we can find $t,t'\in I\setminus\S$ satisfying:
\begin{itemize}
\item  $t/E$ and $t'/E$ are both singletons;
\item  $t,t'$ satisfy the same $\S$-cut, i.e., for each $s\in\S$, $s<t$ iff $s<t'$;
\item $t<w<t'$;
\item $t\in J$,  but $t'\not\in J$.
\end{itemize}
We will produce an automorphism $\pi$ of $(I,<,E,P)$ that fixes $\S$
pointwise and $\pi(t)=t'$. This suffices, as necessarily
$\pi(w)\not\in J$ for any such $\pi$. To produce such a $\pi$, first
choose a suitable proper initial segment $K\subseteq I$ containing
$\S\cup\{t,t'\}$. Note that $K$ is countable, and is a union of
$E$-classes.  Consider the structure $(K/E,<,P)$ formed from the
quotient $K/E$, where $<$ is the inherited linear order and $P(r/E)$
if and only if $P(r)$ held in $(I,<,E,P)$. Now $Th(K/E,<,P)$ is known
to be $\aleph_0$-categorical and eliminate quantifiers. [The theory
is axiomatized by asserting that $<$ is dense linear order with a
least element but no greatest element, and $P$ is a dense/codense
subset.] Thus, there is an automorphism $\pi_0$ of $(K/E,<,P)$ fixing
$\S/E$ pointwise and $\pi(t/E)=t'/E$. As every $E$-class of $K$ is
either a singleton or a countable, dense linear order,
there is
an automorphism $\pi_1$ of $(K,<,E,P)$ fixing $\S$ pointwise and
$\pi_1(t)=t'$ and such that $\pi_1(x)/E = \pi_0(x/E)$. Now the
automorphism $\pi$ of $(I,<,E,P)$ defined by $\pi(u)=\pi_1(u)$ if
$u\in K$, and $\pi(u)=u$ for each $u\in I\setminus K$ is as desired.
$\qed_{\ref{seamless}}$

\medskip

The following construction  codes  a nearly arbitrary subset
$S\subseteq\omega_1$ into an $I^S\in\II^*$. We construct orderings
that avoid the third case of Definition~\ref{suitdef}.

\begin{Construction}\label{IS}  Let $S\subseteq\omega_1$ with $0\not\in S$.  There is
$I^S=(I^S,<,P,E)\in\II^*$ that has a continuous, increasing sequence
$\<J_\alpha:\alpha\in\omega_1\>$ of proper initial segments such that:
\begin{enumerate}
\item If $\alpha\in S$, then $I^S\setminus J_\alpha$ has a minimum element $a_\alpha$ satisfying $P(a_\alpha)$; and
\item  If $\alpha\not\in S$ and $\alpha>0$, then $J_\alpha$ is seamless.
\end{enumerate}
\end{Construction}

\bp  Let $\Ltau=\{<,P,E\}$ and $\A$ be the $\Ltau$-structure with
universe singleton $\{a\}$ with both $P(a)$ and $E(a,a)$ holding. Let
$\B=(\QQ,<,P,E)$, where $(\QQ,<)$ is a countable dense linear order
with no endpoints, $P$ fails everywhere, and all elements are
$E$-equivalent. Combine these to get a (countable) $\Ltau$-structure
$\CC$ formed by the dense/codense (with no endpoints) concatenation
of countably many copies of both $\A$ and $\B$. Finally, take $\D$ to
be the concatenation $\A\conc\CC$.

Using these $\Ltau$-structures as building blocks, form a continuous
sequence of $\Ltau$-structures $J_\alpha$, where $J_\alpha$ is an
$\Ltau$-substructure and an initial segment of $J_\beta$ whenever
$\alpha<\beta$ by:  $J_0$ is the one-element structure $\{\min(I)\}$
with $\neg P(\min(I))$. For $\alpha<\omega_1$ a non-zero limit
ordinal, take  $J_\alpha$ to be the increasing union of
$\<J_\beta:\beta<\alpha\>$. Given $J_\alpha$, form $J_{\alpha+1}$ by
$$J_{\alpha+1}=
 \left\{\begin{array}{ll}

J_\alpha\conc \D &\mbox{if $\alpha\in S$}\\
J_\alpha\conc \CC &\mbox{if $\alpha\not\in S$}
\end{array}
\right. $$ Finally, take $I^S$ to be the increasing union of
$\<J_\alpha:\alpha<\omega_1\>$. $\qed_{\ref{IS}}$



\subsection{Striated models and forcing}

In this section we introduce the notion of a striation of a model - a
decomposition of a model $N$ of $T$ into uncountably many countable
pieces satisfying certain constraints on $\pcl$.  We will show
later how to code stationary sets by specially constructed (forced)
striated models.

\subsubsection{Striated Models}
Fix an atomic  $N\in \At$ and some $I=(I,<,E,P)\in\II^*$.

\begin{Definition}\label{striate}  {\em We say $N$ is {\em striated by $I$} if there are $\omega$-sequences $\<\abar_t:t\in I\>$ satisfying:
\begin{itemize}
\item  $N=\bigcup\{\abar_t:t\in I\}$;  (As notation, for $t\in I$, $N_{<t}=\bigcup\{\abar_j:j<t\}$.)
\item  If $t=\min(I)$, then $\abar_t\subseteq\pcl(\emptyset, N)$;
\item  For $t>\min(I)$, $a_{t,0}\not\in\pcl(N_{<t},N)$;
\item  For each $t$ and $n\in\omega$, $a_{t,n}\in\pcl(N_{<t}\cup\{a_{t,0}\},N)$.
\end{itemize}
}

\end{Definition}

Note:  In the definition above, we allow $a_{s,m}=a_{t,n}$ in some
cases when $(s,m)\neq(t,n)$.  However, if $s<t$, then the element
$a_{t,0}\neq a_{s,m}$ for any $m$.  Also, if
$\pcl(\emptyset,N)=\emptyset$, we do not define $\abar_{\min(I)}$.
Although $E$ and $P$ don't appear explicitly in either
Definition~\ref{striate} or Definition~\ref{catchetal}, $E$ is needed
for the following notations and $P$ plays a major role later.

The idea of our forcing will be to force the existence of a striated
atomic model $N_I$ indexed by a linear order $I\in\II^*$ with
universe $X=\{x_{t,n}:t\in I,n\in\omega\}$.  Such an $N_I$ will have
a `built in' continuous sequence $\<N_\alpha:\alpha\in\omega_1\>$ of
countable, elementary substructures, where the universe of $N_\alpha$
will be $X_\alpha=\{x_{t,n}:t\in J_\alpha, n\in\omega\}$ for some
initial segment $J_\alpha$ of $I$. We start with the assumption that
pseudo-minimal types are not dense so some formula $\delta(x,\fbar)$
has `no pseudo-minimal extension'. We absorb the constants $\fbar$
into the language and use
 the assumption of `no pseudo-minimal
extension'  to make the set
$$\{\alpha\in\omega_1:I\setminus J_\alpha\ \hbox{has a least element\}}$$
(infinitarily) definable.  To make this precise, we introduce some notation.

Suppose that $(I,<,P,E)\in\II^*$ and $N=\{a_{t,n}:t\in I,n\in\omega\}$
is striated by $I$. For any suitable $J\subseteq I$, let $N_J$ denote
the substructure with universe $\{a_{t,n}:t\in J,n\in\omega\}$.
Abusing notation slightly, given any $s\in I\setminus\{\min(I)\}$,
let $$J_{<s}=\{s'\in I:s'<s\ \hbox{and}\ \neg E(s',s)\}$$ Thus,
$J_{<s}$ is a suitable proper initial segment of $I$, and we denote
its associated $L$-structure, $\{a_{t,n}: t\in J_{<s}, n<\omega\}$,
by $N_{<s}$. With this notation, we now describe three relationships
between an element and a substructure of this sort.

\begin{Definition}\label{catchetal}  {\em  Suppose $N$ is striated by $(I,<,P,E)$, $J\subseteq I$ suitable, and $b\in N\setminus N_J$.
\begin{itemize}
\item  {\em $b$ catches $N_J$} if, for every $e\in N$, $e\in\pcl(N_{J}\cup\{b\},N)\setminus N_{J}$ implies $b\in\pcl(N_{J}\cup\{e\},N)$.

\item  {\em $b$ has unbounded reach in $N_J$} if there exists
    $s^*\in J$ such that, letting $A$ denote
 $\pcl(N_{<s^*}\cup\{b\},N) \cap N_J$, for every $s\in J$ with
 $s>s^*$ there is a $c \in A -N_{<s}$.

%
\item  {\em $b$ has bounded effect in $N_J$} if there exists $s^*\in J$ such that $\pcl(N_{<s}\cup\{b\},N)\cap N_J=N_{<s}$ for every $s>s^*$ with $s\in J$.
\end{itemize}
}
\end{Definition}

 Clearly, an element $b$ cannot have both unbounded reach and bounded effect in $N_J$, but the properties are not
 complementary.


\begin{Definition}\label{full}  {\em A model $M$ with uncountable cardinality is said to be {\em full}
if for every $\abar \in M$  every non-algebraic $p \in S_{at}(\abar)$
is realized  $|M|$-times in $M$.
}
\end{Definition}

 The remainder of this section is devoted to the proof of the following Theorem.

\begin{Theorem}  \label{forcing}  Suppose $\delta(x)$ is a complete, non-pseudo algebraic formula with no pseudo-minimal extension.
For every $(I,<,P,E)\in\II^*$ there is a c.c.c.\ forcing $\QQ_I$ such that in $V[G]$, there is a full, atomic $N_I\models T$
striated by $(I,<)$ such that:
\begin{enumerate}
\item  For every suitable initial segment $J\subseteq I$, $N_J\preceq N_I$;
\item  If  $t\in I$ and $P(t)$ holds,
then $a_{t,0}$ catches and has unbounded reach in $N_{<t}$;
\item  If $J\subseteq I$ is seamless, then for every $b\in N_I\setminus N_{J}$,
if $b$ catches $N_J$, then $b$ has bounded effect in $N_J$.
\end{enumerate}
\end{Theorem}

\bp
The hypothesis that $\delta(x)$ has no pseudo-minimal extension means
for every $\phi(x,\abar)$ which implies $\delta(x)$ and is not
pseudoalgebraic there do not exist $\abar^*, c, b$ satisfying the
Definition~\ref{psmdef} of pseudominimality. Replacing $\abar^*, c,
b$ by $\bbar,e,c$, our hypothesis on  $\delta(x)$ translates into the
following statement:

\begin{Fact}  \label{good}  Assume  $\delta(x)$ has no pseudo-minimal extension. For any $M\in\At$, for any $\abar$ from $M$ and any $c\in\delta(M)$ for which $c\not\in\pcl(\abar,M)$,
there are $\bbar$ and $e$ from $M$ such that
\begin{enumerate}
\item  $e\in\pcl(\abar\bbar c,M)\setminus \pcl(\abar\bbar,M)$; but
\item  $c\not\in\pcl(\abar\bbar e, M)$.
\end{enumerate}
\end{Fact}

Fix, for the whole of the proof, some $(I,<,E,P)\in\II^*$.  We wish
to construct an atomic model $N_I\models T$, whose complete diagram
contains variables $\{x_{t,n}:t\in I,n\in\omega\}$, that is striated
by $(I,<)$, and includes $\delta(x_{t,0})$, whenever $I\models P(t)$.
We begin by defining a forcing notion $\QQ_I$ and prove that it
satisfies the c.c.c. Then, we exhibit several collections of subsets
of $\QQ_I$ and prove that each is dense and open. Fact~\ref{good}
will only be used in showing the sets witnessing `unbounded reach'
(i.e., Group F of the constraints) are dense. Finally in
Section~\ref{finarg}, we argue that if $G\subseteq\QQ_I$ is a generic
filter meeting each of these dense open sets, then $V[G]$ will
contain an atomic model $N_I$ of $T$ satisfying the conclusions of
Theorem~\ref{forcing}.

\subsection{The forcing}\label{forcedetail}
Our forcing $\QQ_I$ consists of `finite approximations' of this
complete diagram. The conditions will be complete types in variable
with a specific kind of indexing that we now describe.

\begin{Notation} \label{finiteseq} {\em
A {\em finite sequence $\xbar$ from $\<x_{t,n}:t\in I, n\in\omega\>$}
is {\em indexed by $u$} if it has the form $\xbar=\<x_{t,m}:t\in u,
m<n_t\>$, where $u\subseteq I$ is finite and $1\le n_t<\omega$ for
every $t\in u$.

 Given a finite sequence $\xbar$ indexed by $u$ and $\<n_t:t\in u\>$
and given a proper initial segment $J\subseteq I$, let $u\mr J=u\cap J$ and $\xbar\mr J=\<x_{t,m}:t\in u\mr J, m<n_t\>$.

As well, if $p(\xbar)$ is a complete type in the variables $\xbar$,
then $p\mr J$ denotes the restriction of $p$ to $\xbar\mr J$, which
is necessarily a complete type. For $s\in I$, the  symbols $u\mr{<s}$
and $\xbar\mr {\le s}$ are defined analogously, setting $J=I\mr {<s}$
and $I\mr {\le s}$, respectively.
 If $\xbar$ arises from a type $p$ that we are keeping track of, we write $n_{p,t}$ for $n_t$.
 These various notations may be
combined to yield, for example, $p\mr {\leq s/E}$.}
\end{Notation}


The forcing $\QQ_I$ will consist of finite approximations of a complete diagram of an $L$-structure in the variables
$\{x_{t,\ell}:t\in I,\ell\in\omega\}$.  Recall that
the property, `$a\in\pcl(\bbar)$' is
    enforced by a first order formula; this justifies `say' in the next
     definition.       \medskip
\begin{Definition}[ $(\QQ_{I},\le_\QQ)$] \label{forcedef}  {\em
 $p\in \QQ_I$ if and only if the following conditions hold:
 \begin{enumerate}

\item $p$ is a complete (principal) type with respect to $T$ in the variables $\xbar_p$, which are a finite sequence indexed by $u_p$ and $n_{p,t}$
(when $p$ is understood we sometimes write $n_t$);

\item If  $t \in u_p$ and $P(t)$ holds, then $p\vdash\delta(x_{t,0})$;

\item If $t=\min(I)$, then $p$ `says' $\{x_{t,n}:n<n_t\}\subseteq\pcl(\emptyset)$;

\item  If $p$ `says' $x_{t,0}\in\pcl(\emptyset)$, then $t=\min(I)$;

\item  For all $t\in u_p$, $t\neq\min(I)$, $p$ `says' $x_{t,0}\not\in\pcl(\xbar_p\mr {<t})$; and
\item  For all $t\in u_p$ and $m<n_t$, $p$ `says' $x_{t,m}\in\pcl(\xbar_p\mr {<t}\cup\{x_{t,0}\})$.

\end{enumerate}

For $p,q\in\QQ_I$, we define
$p\le_{\QQ_I} q$ if and only if $\xbar_p\subseteq\xbar_q$ and the complete type $p(\xbar_p)$ is
 the restriction of $q(\xbar_q)$ to $\xbar_p$.}
\end{Definition}
\medskip

We begin with some easy observations.

\begin{Lemma}  \label{truncate}
For every $p\in\QQ_I$ and every proper initial segment $J\subseteq
I$, $p\mr J\in\QQ_I$ and $p\mr J\le_{\QQ_I} p$.
\end{Lemma}

\begin{Lemma}  \label{extendauto}
Every automorphism $\pi$ of $(I,<,E,P)$ naturally extends to an automorphism $\pi'$ of $\QQ_I$
via the mapping $x_{t,n}\mapsto x_{\pi(t),n}$.
\end{Lemma}


\begin{Lemma}  \label{chaincondition} Suppose $p\in\QQ_I$ and $u_p\neq\emptyset$.  Enumerate $u_p=\{s_i:i<d\}$ with $s_i<_I s_{i+1}$ for each $i$.
For any $M\in\At$ and any $\bbar$ from $M$ realizing $p(\xbar_p)$,
there is a sequence $M_0\preceq M_1\preceq\dots\preceq M_{d-1}=M$ of elementary substructures of $M$ satisfying:
\begin{itemize}
\item  For each $i<d$, $\bbar\mr{< s_i}\subseteq M_i$; and
\item  For $0<i<d$, $b_{s_i,0}\in M_i\setminus M_{i-1}$.
\end{itemize}
\end{Lemma}

\bp  By induction on $d=|u_p|$.  For $d=0,1$ there is nothing to
prove, so assume $d\ge 2$ and the Lemma holds for $d-1$. Fix any
$M\in\At$ and choose any realization $\bbar$ of
$p(\xbar_p)$ in $M$.  Clearly, the subsequence
$\abar:=\bbar\mr{<s_{d-1}}$ realizes the restriction
$q:=p\mr{<s_{d-1}}$.  As $b_{s_{d-1},0}\not\in\pcl(\abar,M)$, there
is $M_{d-2}\preceq M$ such that $\abar$ is from $M_{d-2}$, but
$b_{s_{d-1},0}\not\in M_{d-2}$.  Then complete the chain by applying
the inductive hypothesis to $M_{d-2}$ and $q$.
$\qed_{\ref{chaincondition}}$

\medskip

The `moreover' in the following lemma emphasizes that in proving
density we are showing how to assign levels to a elements of a finite
sequence in a model which need not be striated.

\begin{Lemma}  \label{consistent} Suppose $J\subseteq I$ is an initial segment and $p,q\in\QQ_I$ satisfy $p\mr J\le_\QQ q$ and $u_q\subseteq J$.
Then there is $r\in\QQ_I$ with $\xbar_r=\xbar_p\cup\xbar_q$,
$r\ge_\QQ p$ and $r\ge_\QQ q$.  Moreover, if $M\in\At$, $\abar$
realizes $p\mr{J}$, $\abar\bbar$ realizes $p$, and $\abar\cbar$
realizes $q$, then $\abar\bbar\cbar$ realizes $r$.
\end{Lemma}

\bp  If $u_p=\emptyset$, then take $r=q$, so assume otherwise.
Choose any $M\in\At$ and fix a realization $\bbar$ of $p(\xbar_p)$ in $M$.  Let $\abar=\bbar\mr J$.
Write $u_p=\{s_i:i<d\}$ with $s_i<_I s_{i+1}$ for each $i$.
 Apply Lemma~\ref{chaincondition} to $M$ and $\bbar$ and choose $\ell<d$ least such that $\abar\subseteq M_\ell$.
 As $q(\xbar_q)$ is generated by a complete formula and $\abar\subseteq M_\ell$, there is
$\cbar\subseteq M_\ell$ such that $\abar\cbar$ (when properly
indexed) realizes $q$. Now define $r(\xbar_r)$ to be the complete
type of $\bbar\cbar = \abar\bbar\cbar$ in $M$ in the variables
$\xbar_r=\xbar_p\cup\xbar_q$. $\qed_{\ref{consistent}}$

\begin{Claim}\label{ccc}
$(\QQ_I,\le_\Q)$ has the c.c.c.
\end{Claim}

Proof.  Let $\{ p_i\mcolon i< \aleph_1\}\subseteq\QQ_I$ be  a collection of
conditions.  We will find $i\neq j$ for which $p_i$ and $p_j$ are compatible.
We successively reduce this set maintaining its uncountability. By
the $\Delta$-system lemma we may assume that there is a single $u^*$
such that for all $i,j$, $u_{p_i} \cap u_{p_j} = u^*$. Further, by the
pigeonhole principle we can assume that for each $t\in u^*$, $n_{p_i,t}
= n_{p_j,t}$.   We can use pigeon-hole again to guarantee that all
the $p_i$ and $p_j$ agree on the finite set of shared variables. And
finally, since $I$ is $\aleph_1$-like we can choose an uncountable
set $X$ of conditions such that for $i<j$ and $p_i,p_j \in X$ all
elements of $u^*$ precede anything in any $u_{p_i}\setminus u^*$ or $u_{p_j}\setminus u^*$
and that all elements of $u_{p_i}\setminus u^*$ are less that all elements of
$u_{p_j} \setminus u^*$.

Finally, choose any $i<j$ from $X$.  Let $J=\{s\in I:s\le\max(u_{p_i})\}$.
By Lemma~\ref{consistent} applied to $p_i$ and $p_j$ for this choice of $J$, we conclude that $p_i$ and $p_j$ are compatible.
$\qed_{\ref{ccc}}$

\bigskip

Recall that a set $X \subseteq \QQ_I$ is dense if for every $p\in \QQ_I$
there is a $q \in X$ with $q \geq p$ and $X \subseteq \QQ_I$ is open if
for every $p\in X$ and $q \geq p$, then
 $q \in X$.

In the remainder of Section~\ref{forcedetail}  we list the crucial
`constraints', which are sets of conditions,  and we prove each of
them to be dense and open in $\QQ_I$.


\medskip\par\noindent
{\bf A.  Surjectivity}
Our first group of constraints ensure that for any generic $G\subseteq\QQ_I$,
for every $(t,n)\in I\times\omega$, there is $p\in G$ such that $x_{t,n}\in \xbar_p$.
To enforce this, for
any $(t,n)\in I\times\omega$, let
$$\A_{t,n}=\{p\in\QQ_I: x_{t,n}\in\xbar_p\}$$

%

\begin{Claim}  \label{surjective}
\begin{enumerate}
\item  For every $t\in I\setminus\{\min(I)\}$ and every $n\in\omega$, $\A_{t,n}$ is dense and open;
\item  If $\pcl(\emptyset)\neq\emptyset$, then $\A_{\min(I),n}$ is dense and open for every $n\in\omega$.
\end{enumerate}
Moreover, in either case, given $(t,n)\in I\times\omega$ and any $p\in\QQ_I$, there is $q\in\A_{t,n}$ with $q\ge_\QQ p$
and $u_q=u_p\cup\{t\}$.
\end{Claim}


\bp  Each of these sets are trivially open.  We first establish
density for (1) and (2) when $n=0$.  For $t = \min(I)$, (1) is
vacuous. For (2), choose any $p\in\QQ_I$. If
$x_{\min(I),0}\in\xbar_p$, there is nothing to prove, so assume it is
not. Pick any $M\in\At$. Choose $\bbar$ from $M$ realizing $p$ and
choose $a\in\pcl(\emptyset,M)$. Then define $q$ by
$\xbar_q=\xbar_p\cup\{x_{\min(I),0}\}$ and $q(\xbar_q)=\tp(\bbar
a,M)$. Next, we show that $\A_{t,0}$ is dense for every $t>\min(I)$.
To see this, choose any $p\in\QQ_I$.  If $t\in u_p$, then necessarily
$x_{t,0}\in\xbar_p$, so there is nothing to prove.  Thus, assume
$t\not\in u_p$. Take $J=\{s\in I:s< t\}$.  Pick $M\in\At$ and choose
a realization $\abar$ of $p\mr J$ in $M$.

%
%
%
%

 As $\delta$ is not pseudo-algebraic, by Lemma~\ref{outside} there is $b\in M$ realizing
$\delta$ with $b\not\in\pcl(\abar,M)$.  Let $q\in\QQ_I$ be defined by
$\xbar_q=\xbar_p\mr J\cup\{x_{t,0}\}$ and the complete type
$q(\xbar_q)=\tp(\abar b,M)$.  Then $q\ge_\QQ p\mr J$ and by
Lemma~\ref{consistent}, there is $r\in\QQ$ with $r\ge_\QQ q$ and
$r\ge_\QQ p$.  Visibly, $r\in\A_{t,0}$.

Next, we prove by induction on $n$ that if $\A_{t,n}$ is dense, then so is $\A_{t,n+1}$.  But this is trivial.  Fix $t$ and choose
$p\in\QQ_I$ arbitrarily.  By our inductive hypothesis, there is $q\ge p$ with $x_{t,n}\in \xbar_q$.  If $x_{t,n+1}\in\xbar_q$,
there is nothing to
prove, so assume otherwise.  Then, necessarily, $n_{q,t}=n+1$.  Let $r$ be the extension of
 $q$ with  $\xbar_r=\xbar_q\cup\{x_{t,n+1}\}$
and $r(\xbar_r)$ the complete type generated by $q(\xbar_q)\cup\{x_{t,n+1}=x_{t,n}\}$.

The final sentence holds by inspection of the proof above.
$\qed_{\ref{surjective}}$

\medskip\par\noindent {\bf B.  Henkin witnesses}

For every $t\in I$, for every finite sequence $\xbar$ (indexed as in
Notation~\ref{finiteseq}) from $I\mr{<t}\times\omega$, and for
every $L$-formula $\phi(y,\xbar)$, $\B_{\phi,t}$ is the set of $p \in
\QQ$ such that:
\begin{enumerate}
\item $\xbar\subseteq\xbar_p$; and
    \item Some $s\in u_p$ and   $m<n_{p,s}$ satisfy $s<t$ and  $p(\xbar_p) \vdash (\exists
        y)\phi(y,\xbar) \rightarrow
        \phi(x_{s,m},\xbar).$

 \end{enumerate}

 \begin{Claim}  \label{Henkin}  For each $t\in I$, finite sequence $\xbar$ from $I\mr{<t}\times\omega$, and $\phi(y,\xbar)$, $\B_{\phi,t}$ is dense and open.
\end{Claim}


\bp   Fix $t\in I$ and $\phi(y,\xbar)$ as above.  Choose any $p\in\QQ_I$.  By using Claim~\ref{surjective} and extending $p$ as needed, we may assume $\xbar\subseteq
\xbar_p$.  Let $q$ denote $p\mr{<t}$.  Then $q\in\QQ_I$ and $q\le_\QQ p$ by Lemma~\ref{truncate}.  As $\xbar\subseteq I_{<t}\times \omega$, $\xbar\subseteq\xbar_q$,
so by adding dummy variables to $\phi$ we may assume $\xbar=\xbar_q$.
 Choose any $M\in\At$ and any realization $\bbar$ of $q$.
 There are now a number of cases.

 Case 1:  $M\models\neg\exists y\phi(y,\bbar)$.  Then as $q(\xbar)$ generates  a complete type, $q\vdash\neg\exists y \phi(y,\xbar_q)$, hence $p\in\B_{\phi,t}$.

 So, we assume this is not the case.  Fix a witness $c\in M$ such that $M\models\phi(c,\bbar)$.
There are now several cases depending on the complexity of $c$ over $\bbar$.
In each of them, we will produce $r\ge_\QQ q$ with $u_r\subseteq I\mr{<t}$ and $r(\xbar_r)\vdash\exists y\phi(y,\xbar)$.

Case 2:  $c\in\pcl(\emptyset,M)$.  If $\min(I)\not\in q$, then let $\xbar_r=\xbar_q\cup\{x_{\min(I),0}\}$ and if $\min(I)\in q$,
then let $\xbar_r=\xbar_q\cup\{x_{\min(I),m}\}$, where $m=n_{q,\min(I)}$.  Regardless,  put $r(\xbar_r)=\tp(\bbar c,M)$.

Case 3:  $c\not\in\pcl(\bbar,M)$.  Choose $s^*>u_q$ with $s^*<t$.
  Let $\xbar_r=\xbar_q\cup\{x_{s^*,0}\}$ and again take $r(\xbar_r)=\tp(\bbar c,M)$.  It is easily checked that $r\in\QQ_I$.

Case 4:  $c\in\pcl(\bbar,M)\setminus\pcl(\emptyset,M)$.  For each $s\in u_q$, let $\xbar\mr{\le s}$ be the subsequence of $\xbar$ consisting of all $x_{t,m}\in\xbar$ with
$t\le s$, and let $\bbar\mr{\le s}$ be the corresponding subsequence of $\bbar$.   Using this as notation, choose $t^*\in u_q\setminus\{\min(I)\}$ least such that
$c\in\pcl(\bbar\mr{\le t^*},M)$.
Again, let $\xbar_r=\xbar_q\cup\{x_{t^*,m}\}$, where $m=n_{q,t^*}$, and
let $r(\xbar_r)=\tp(\bbar c,M)$.  As in the case above, it is easily verified that $r\in\QQ_I$.

Now, in any of Cases 2,3,4, by Lemma~\ref{consistent} we can find $p^*\ge_{\QQ} p$ and $p^*\ge_{\QQ} r$.
$\qed_{\ref{Henkin}}$

 \medskip\par\noindent {\bf C.  Fullness}
Suppose $\xbar$ is a finite sequence (indexed as in
Notation~\ref{finiteseq}), $t\in I$, and $\phi(y,\xbar)$ is an
$L$-formula such that $\phi(y,\xbar)$ `says' `$y$ is not
pseudo-algebraic over $\xbar$.'

$$\CC_{\phi,t}=\{p\in\QQ_I:\ \hbox{there is}\ s>t, s\in u_p,\xbar\subseteq\xbar_p, p\vdash\phi(x_{s,0},\xbar)\}$$

\begin{Claim}  \label{fullness}  Each is $\CC_{\phi,t}$ is dense and open.
\end{Claim}

\bp  Fix $\phi(y,\xbar)$ and $t$, and choose any $p\in\QQ_I$.  By
extending $p$ as needed, by Claim~\ref{surjective} we may assume
$\xbar\subseteq\xbar_p$. Choose any countable $M\in\At$ and choose
any realization $\bbar$ of $p(\xbar_p)$ in $M$.  As $\phi(y,\bbar)$
is not pseudo-algebraic, there is $N\in\At$, $N\succeq M$, and $c\in
N\setminus M$ satisfying $N\models\phi(c,\bbar)$.  Choose any $s\in
I$ such that $s>\max(u_p)$ and $s>t$ with $I\models\neg P(s)$. Define
$q$ by:  $\xbar_q=\{x_{s,0}\}\cup\xbar_p$ and
$q(\xbar_q)=\tp(c\bbar,N)$.  Then $q\ge_\QQ p$ and $q\in \CC_{\phi,t}$.
$\qed_{\ref{fullness}}$

%

\medskip\par\noindent {\bf D+E.  Determining level}
         The definition of the forcing implies that $x_{t,n}$ is pseudo-algebraic
         over $\xbar_p\mr{<t}\cup\{x_{t,0}\}$ for any $p\in\QQ_I$ with $x_{t,n}\in\xbar_p$,
         but it might also be algebraic over some smaller finite sequence (at a lower level).  If this occurs,
          we `adjust the level' by finding some  $s<t$ and $m$  and insisting that
         $x_{t,n}=x_{s,m}$.  To make this precise involves defining two families of constraints
          and showing that each is dense and open. The first family
          is actually a union of two.

$\D_{t,n}= \D^1_{t,n}\cup \D^2_{t,n} $ where
\begin{enumerate}
\item  $\D^1_{t,n}=\{p:x_{t,n}\in\xbar_p$ and $p$ `says'
    $x_{t,0}\in\pcl(\xbar_p\mr{<t}\cup\{x_{t,n}\})\}$;

 \item  $\D^2_{t,n}=\{p:x_{t,n}\in\xbar_p$, there are $s\in u_p$,
     $s<t$, and $m<n_{p,s}$ such that $p(\xbar_p)\vdash
     x_{t,n}=x_{s,m}\}$.

    \end{enumerate}

The second family is parameterized by $\xbar, t, n$. Let $\xbar$ be
any finite sequence (cf.\ Notation~\ref{finiteseq}) indexed by $u$
with $s=\max(u)<t$.
$$\E_{t,n,\xbar}=\{p\in\QQ_I:\xbar\cup\{x_{t,n}\}\subseteq\xbar_p\ \hbox{and {\bf either} $p$ `says' $x_{t,n}\not\in\pcl(\xbar)$
{\bf or} $p$ `says' $x_{t,n}=x_{s,m}$ for some $m$}\}$$

\begin{Claim} \label{level1} For all $(t,n)\in I\times\omega$ and for all finite sequences $\xbar$ indexed by $u$ with $\max(u)<t$, $\E_{t,n,\xbar}$ is dense and open.
\end{Claim}

\bp  Once more, `Open' is clear.  Let $s=\max(u)$.  Given any
$p\in\QQ_I$, by iterating Claim~\ref{surjective} we may assume
$\xbar\cup\{x_{t,n}\}\subseteq\xbar_p$.  If $p$ `says'
$x_{t,n}\not\in\pcl(\xbar)$, then $p\in\E_{t,n,\xbar}$, so assume $p$
`says' $x_{t,n}\in\pcl(\xbar)$.  From our conditions on $\xbar$, this
implies $x_{t,n}\in\pcl(\xbar_p\mr{\le s})$.  So put $m=n_{p,s}$, let
$\xbar_q=\xbar_p\cup\{x_{s,m}\}$ and let $q(\xbar_q)$ be the complete
type generated by $p(\xbar_p)\cup \hbox{`$x_{t,n}=x_{s,m}$'}$.
$\qed_{\ref{level1}}$

    \begin{Claim}  \label{level}

For every $t\in I\setminus\{\min(I)\}$ and every $n\in\omega$, $\D_{t,n}$ is dense and open.
\end{Claim}

\bp  Choose any $p\in\QQ_I$.  By Claim~\ref{surjective} we may assume $x_{t,n}\in\xbar_p$.
Choose any $M\in\At$ and choose $\bbar$ in $M$ realizing $p$.    There are now several cases.

Case 1.  If $b_{t,0}\in\pcl(\bbar\mr{<t}\cup\{b_{t,n}\})$, then $p\in\D^1_{t,n}$, so assume this is not the case.

Case 2.  If $b_{t,n}\in\pcl(\emptyset,M)$ and $\min(I)\not\in u_p$, then define $q$ by $\xbar_q=\xbar_p\cup\{x_{\min(I),0}\}$
and $q(\xbar_q)=\tp(\bbar b_{t,n},M)$.

Case 3.  If $b_{t,n}\in\pcl(\bbar_{\le s},M)$ for some $s\in u_p$, $s<t$, then define $q$ by $\xbar_q=\xbar_p\cup\{x_{s,m}\}$ (where $m=n_{p,s}$)
and $q(\xbar_q)$ be the extension of $p(\xbar_p)$ by `$x_{t,n}=x_{s,m}$.'

%

Case 4. If none of the previous cases occur, choose  $s^*<t$ with
$s^*> u_p \cap I_{<t}$, $I\models\neg P(s^*)$. Define $q$ by
$\xbar_q=\xbar_p\cup\{x_{s^*,0}\}$ and $q(\xbar_q)=\tp(\bbar
b_{t,n},M)$ (i.e. $x_{s^*,0} = x_{t,n}$).  Now since Case~1 fails,
$q$ satisfies Condition 5) in the definition of $\QQ_I$ at level $t$,
and since Case 3 fails, Condition~5) holds at level $s^*$. And in q,
Condition~6) holds for $x_{t,n}$ since $b_{t,n} =b_{s^*,0}$.  The
other conditions are inherited from $p$, so $q\in\QQ_I$.
 $\qed_{\ref{level}}$

   \medskip\par\noindent {\bf F.  Achieving unbounded reach}


Suppose $s_0/E<s_1<t$ are from $I$ with $I\models P(t)$,
$s_0\neq\min(I)$, and $I\models\neg P(s_0)$ (so $s_0/E$ is infinite
and dense).

$\F_{t,s_0,s_1}$ is the set of $p \in \QQ_I$ such that there exists
$s_2 \in u_p$ with $s_1<s_2 <t$ such that (recalling
Notation~\ref{finiteseq}) $p$ `says'
$$x_{s_2,0} \in \pcl(\{x_{t,0}\} \cup \xbar_p\mr{\le s_0/E}).$$



\begin{Claim}  \label{blocking}  Each $\F_{t,s_0,s_1}$ is dense and open.
\end{Claim}

\bp  Open is clear.  Choose any $p\in\QQ_I$.  By Claim~\ref{surjective} we may assume $x_{t,0}\in\xbar_p$.
 By Lemma~\ref{truncate} we have the sequence of extensions:
$$p\mr{\le s_0/E}\ \le_\QQ \ p\mr{<t}\ \le_\QQ\  p\mr{\le t}\ \le_\QQ \ p.$$
Fix $M\in\At$ and choose sequences $\abar$, $\dbar$, $\cbar$ from $M$
such that $\abar\dbar\cbar$ realizes $p\mr{\le t}$, with $\abar$
realizing $p\mr{\le s_0/E}$ and $\cbar$ realizing $p\mr{=t}$.  Let
$c_0\in\cbar$ be the interpretation of $x_{t,0}$.  Thus,
$M\models\delta(c_0)$ and $c_0\not\in\pcl(\abar,M)$. Using
Fact~\ref{good}, choose $\bbar$ and $e$ from $M$ such that
$e\in\pcl(\abar\bbar c_0,M)\setminus\pcl(\abar\bbar,M)$, but
$c_0\not\in\pcl(\abar\bbar e,M)$. We will find conditions in $\QQ$
that assign levels to $\bbar$ and $e$ to satisfy $\F_{t,s_0,s_1}$.

As the class $s_0/E$ has no last element, by using Claim~\ref{Henkin}
(Henkin witnesses) $\lg(\bbar)$ times, we can construct $q\in\QQ_I$,
$q\ge_\QQ p\mr{\le s_0/E}$ satisfying $q(\xbar_q)=\tp(\abar\bbar,M)$
and $u_q\subseteq I\mr{\le s_0/E}$.

Next, by Lemma~\ref{consistent} there is $q_1\ge_\QQ q$, $q_1\ge_\QQ
p\mr{<t}$, and $u_{q_1}\subseteq I\mr{<t}$. By Lemma~\ref{consistent}
again, there is $q_2\ge_\QQ q_1$, $q_2\ge_\QQ p\mr{\le t}$, and
$u_{q_2}\subseteq I\mr{\le t}$. Indeed, by the `Moreover' clause of
Lemma~\ref{consistent}, we may additionally assume that
$q_2(\xbar_{q_2})=\tp(\abar\bbar\dbar\cbar,M)$ (and so
$q_1(\xbar_{q_1})=\tp(\abar\bbar\dbar,M)$).

Now, choose $s_2\in I$ such that $I\models\neg P(s_2)$, $s_1<s_2<t$, and $s_2>s$ for every $s\in u_{q_1}$.  Define $r$ by
$\xbar_r=\xbar_{q_2}\cup\{x_{s_2,0}\}$ and $r(\xbar_r)=\tp(\abar\bbar\dbar\cbar e,M)$.     It is easily checked that $r\in\QQ_I$
and visibly, $r\ge_\QQ q_2$.  As well, $r\in\F_{t,s_0,s_1}$.

Finally,  by a final application of Lemma~\ref{consistent}, since
$u_r\subseteq I\mr{\le t}$ and $r\ge_\QQ p\mr{\le t}$, there is
$p^*\ge_\QQ p$ with $p^*\ge_\QQ r$.  As $p^*\in\F_{t,s_0,s_1}$,  we
conclude that $\F_{t,s_0,s_1}$ is dense. $\qed_{\ref{blocking}}$

\medskip\par

\subsection{Proof of Theorem~\ref{forcing}}\label{finarg}

Given a linear order $I$ we construct a model $N=N_I$ of the theory
$T$.  That is,  we verify that the forcing $(\QQ_I,\le_\QQ)$
satisfies the conclusions of Theorem~\ref{forcing}. Suppose
$G\subseteq\QQ_I$ is a filter meeting every dense open subset. Let
$$X[G]=\bigcup\{p(\xbar_p):p\in G\}$$ Because of the dense subsets
$\A_{t,n}$, $X[G]$  describes a complete type in the variables
$\{x_{t,n}:t\in I,n\in\omega\}$.\footnote{If
$\pcl(\emptyset)=\emptyset$, then $X[G]$ is in the variables
$\{x_{t,n}:n\in\omega, t\in I\setminus\{\min(I)\}\}$. For clarity of
exposition, we will assume that $\pcl(\emptyset)\neq\emptyset$.}
Intuitively, we want to build a with domain given by these variables.
 But the Level conditions, Claim~\ref{level} introduced  a natural equivalence relation $\sim_G$ on $X[G]$ defined by
$$x_{t,n}\sim_G x_{s,m}\qquad\hbox{if and only if}\qquad X[G] \ \hbox{`says'}\ x_{t,n}=x_{s,m}$$

Let $N[G]$ be the $\tau$-structure with universe $X[G]/\sim_G$.   Each element of $N[G]$ has the form $[x_{t,n}]$, which is the equivalence class
of $x_{t,n}$ (mod $\sim_G$).  
 As each $p\in\QQ_I$ describes a complete (principal) formula with respect to $T$, $N[G]$ is an atomic set.  As well, it follows from Claim~\ref{Henkin}
that $N[G]\models T$.

For each $t\in I$ such that $P(t)$ holds, let $N_{<t}=\{[x_{w,n}]:$ some $x_{s,m}\in[x_{w,n}]$ with $s<t\}$.
Similarly, for each $s\in I\setminus\{\min(I)\}$ with $\neg P(s)$, let
$N_{<s}=\{[x_{w,n}]:w/E < s/E\}$.

By repeated use of Claim~\ref{Henkin}, both $N_{<t}$ and $N_{<s}$ are
elementary substructures of $N[G]$. Note that  $N_{<s'}=N_{<s}$
whenever $E(s',s)$.

For simplicity, let $a_{w,n}\in N[G]$ denote the class $[x_{w,n}]$.
Given any $(w,n)$, if there is a least $s\in I$ such that
$a_{w,n}=a_{s,m}$ for some $m\in\omega$, then we say {\em $a_{w,n}$
is on level $s$}. For an arbitrary $(w,n)$, a least $s$ need not
exist, but it does in some cases. In particular,
Definition~\ref{forcedef}.5 and the level constraint
($\E_{w,0,\xbar}$) imply that any $a_{w,0}$ is on level $w$ for any
$w\in I$. As well, because of the Level constraints (group $D + E$)
for any $t$ such that $P(t)$ holds and for any $n>0$,
$$\hbox{$a_{t,n}$ is on level $t$ if and only if   $a_{t,0}\in\pcl(N_{<t}\cup\{a_{t,n}\},N[G])$}$$


As $|I|=\aleph_1$ and the fact that each $a_{t,0}\not\in\pcl(N_{<t}, N[G])$, $||N[G]||=\aleph_1$.
Finally, it follows from the density of the `Fullness conditions' that $N[G]$ is full.


It remains to verify that $N[G]$ satisfies the three conditions of Theorem~\ref{forcing}.  First, for any initial segment $J\subseteq I$ without a maximum element (in particular, for any
suitable $J$) the density of the Henkin conditions offered by Claim~\ref{Henkin} and the Tarski-Vaught criterion imply that $N_J\preceq N[G]$.

%
Second, suppose $t\in I$ and $P(t)$ holds.  We show that $a_{t,0}$ catches and has
unbounded reach in $N_{<t}$.  Note that since $I\mr{<t}$ is suitable,
 $N_{<t}\preceq N[G]$, hence $\pcl(N_{<t},N[G])=N_{<t}$.
 To see that $a_{t,0}$ catches $N_{<t}$, choose any
$a_{s,m}\in\pcl(N_{<t}\cup\{a_{t,0}\},N[G])\setminus N_{<t} $.
 By taking an appropriate finite sequence $\xbar$ witnessing the pseudo-algebraicity,
 the density of the constraints $\E_{s,m,\xbar}$ allow us to assume $s\le t$.  However, if $s<t$, then
 we would have $a_{s,m}\in N_{<t}$.  Thus, the only possibility is that $(s,m)=(t,n)$ for some $n\in\omega$
 and that $a_{t,n}$ is on level $t$.  It follows from the displayed remark above that
$a_{t_0}\in\pcl(N_{<t}\cup\{a_{t,n}\},N[G])$.  Thus, $a_{t,0}$ catches $N_{<t}$.
We also argue that $a_{t,0}$ has unbounded reach in $N_{<t}$.
To see this, choose any  $s_0<t$, $s_0\neq\min(I)$ with $I\models\neg P(s_0)$.  For any $s_1$ satisfying $s_0/E<s_1/E<t$,
choose $p\in G\cap\F_{t,s_0,s_1}$ and choose $s_2\in u_p$ from there.  Now,  the element $a_{s_2,0}\in \pcl(N_{<s_0}\cup\{a_{t,0}\},N[G])$.
As well, since $s_1/E<s_2/E<t/E$, $a_{s_2,0}\not\in N_{<s_1}$, so $a_{t,0}$ has unbounded reach in $N_{<t}$.


It remains to verify (3) of Theorem~\ref{forcing}.  Choose a seamless
$J\subseteq I$ and suppose some $b\in N[G]\setminus N_J$ catches
$N_J$. Say $b$ is $a_{t^*,n}$, where necessarily $t^*\in I\setminus
J$. We must show $b$ has bounded effect in $N_J$. By the fundamental
theorem of forcing, there is $p\in G$ such that
$$p\Vdash a_{t^*,n}\ \hbox{catches}\ N_J.$$
Thus, among other things, $p\Vdash$ `$a_{t^*,n}\neq a_{s,m}$' for all
$s\in J$, $m\in\omega$.
%

Choose any $s^*\in J$ such that $s^*>s$ for every $s\in u_p\cap J$.

\begin{Claim}  \label{forceit}
$p\Vdash \pcl(\{b\}\cup N_{<s^*}[\Gdot],N[\Gdot])\cap N_J[\Gdot]\subseteq N_{<s^*}[\Gdot]$.
\end{Claim}

%

\bp  If not, then there is $q\in\QQ_I$ satisfying $q\ge p$ and a finite $A\subseteq N_{<s^*}[\Gdot]$ such that
$$q\Vdash \pcl(Ab,N[\Gdot])\cap N_J[\Gdot]\not\subseteq N_{<s^*}[\Gdot]$$
Without loss, we may assume that if $a_{t,m}\in A$, then $t\in u_q$.
As $J$ is seamless, by Lemma~\ref{seamless}, choose an automorphism
$\pi$ of $(I,<,E,P)$ such that $\pi\mr{\ge\min (u_p\setminus J)}=id$;
$\pi(t^*)=t^*$; $\pi\mr{u_p}=id$; $\pi\mr{u_q\cap I_{<s^*}}=id$, but
$\pi(s^*)\not\in J$. By Lemma~\ref{extendauto}, $\pi$ extends to an
automorphism $\pi'$ of $\QQ_I$ given by $x_{t,m}\mapsto
x_{\pi(t),m}$. By our choice of $\pi$, $\pi'(p)=p$.  While $\pi'(q)$
need not equal $q$, we do have $p\le\pi'(q)$. Now
$$\pi'(q)\Vdash \pcl(Ab,N[\Gdot])\cap N_{\pi(J)}[\Gdot]\not\subseteq
N_{<\pi(s^*)}[\Gdot]$$ But this contradicts $p\Vdash a_{t^*,n}\
\hbox{catches}\ N_J.$  [To see this, choose $H$ generic with
$\pi'(q)\in H$, hence also $p\in H$. Choose $e\in (\pcl(Ab,N[H])\cap
N_{\pi(J)}[H])\setminus N_{<\pi(s^*)}[H]$.  As $A\subseteq N_J[H]$,
$e\in\pcl(N_J[H]\cup\{b\},N[H])$.  Moreover,
as $N_J[H]\preceq
N_{<\pi(s^*)}[H]$, $e\not\in N_J[H]$.  But, since $N_J[H]\cup\{e\}\subseteq N_{\pi(J)}[H]$
and $b\not\in N_{\pi(J)}[H]$, it follows that $b\not\in \pcl(N_J[H]\cup\{e\},N[H])$.
That is, $e$ witnesses that $b$ does not catch $N_J[H]$.]
$\qed_{\ref{forceit}}$


As Claim~\ref{forceit} holds for any sufficiently large $s^*\in J$,  $a_{t,n}$ has bounded effect in $N_J$
This concludes the proof Theorem~\ref{forcing}.$\qed_{\ref{forcing}}$

\section{Proof of Theorem~\ref{big}}  \label{5}


Now we prove the main theorem, Theorem~\ref{big}, by using the
transfer lemma, Theorem~\ref{transfer1} to move from coding a model
by $S$ in $M[G]$ (Theorem~\ref{forcing}) to $2^{\aleph_1}$ models in
$V$.

We  prove Theorem~\ref{big} under the assumption that a countable,
transitive model $(M,\epsilon)$ of a suitable finitely axiomatizable
subtheory of
 ZFC exists.\footnote{Alternatively, one could use the fragment  $ZFC^0$ of \cite{BL}.}
%
As the existence of the latter is
provable from ZFC (using the Reflection Theorem) we obtain a proof of
Theorem~\ref{big} in ZFC.

As the pseudo-minimal types are not dense, we can find a complete
formula $\delta(x,\abar)$ that is not pseudo-algebraic, but has no
pseudo-minimal extension. As having $2^{\aleph_1}$ models is
invariant under naming finitely many constants, we absorb $\abar$
into the signature and write $\delta(x)$ for this complete formula.

Fix a countable, transitive model $(M,\epsilon)$ of ZFC with
$T,\tau\in M$ and we begin working inside it.  In particular, choose
$S\subseteq\omega_1^M\setminus\{0\}$ such that $$(M,\epsilon)\models `\hbox{$S$ is
stationary/costationary'}$$  Next, perform Construction~\ref{IS}
inside $M$ to obtain $I=(I^S,<,P,E)\in\II^*$.

Next, we force with the c.c.c.\ poset $\QQ_{I^S}$ and find
$(M[G],\epsilon)$, where $G$ is a generic subset of $\QQ_{I^S}$. As
the forcing is c.c.c., it follows that all cardinals as well as
stationarity, are preserved, Thus, $\omega_1^{M[G]}=\omega_1^M$ and
$(M[G],\epsilon)\models `\hbox{$S$ is stationary/costationary'}$.


As Construction~\ref{IS} is absolute, $I^{M[G]}=I^{M}=I^S$. According
to Theorem~\ref{forcing}, inside $M[G]$ there is an atomic, full
$N_I\models T$ that is striated according to $(I^S,<,P,E)$. Write the
universe of $N_I$ as $\{a_{t,n}:t\in I^S,n\in\omega\}$. Inside $M[G]$
we have the mapping $\alpha\mapsto J_\alpha$ given by
Construction~\ref{IS}.  For every $\alpha\in\omega_1^{M[G]}$, let
$N_\alpha$ be the $\tau$-substructure of $N_I$ with universe
$\{a_{t,n}:t\in J_\alpha, n\in\omega\}$.  It follows from
Theorem~\ref{forcing} and Construction~\ref{IS} that for every
non-zero $\alpha\in\omega_1^{M[G]}$:
\begin{itemize}
\item  $N_\alpha\preceq N_I$;
\item If $\alpha\in S$, then $I^S\setminus J_\alpha$ has a least element $t(\alpha)$ and $a_{t(\alpha),0}$ both catches and has unbounded reach in $N_\alpha$;

\item If $\alpha\not\in S$, then every $b\in N_I\setminus N_\alpha$ that catches $N_\alpha$ has bounded effect in $N_\alpha$.
\end{itemize}

%
%

%

Now, still working inside $M[G]$, we identify a 3-sorted structure $N^*$ that encodes this information.
  The vocabulary of $N^*$ will be
$$\tau^*=\tau\cup\{U,V,W,<_U,<_V,P,E,R_1,R_2\}.$$
$N^*$ is the $\tau^*$-structure in which
\begin{itemize}
\item $\{U,V,W\}$ are unary predicates that partition the
    universe;
\item  $(U^{N^*},<_U)$ is  $(\omega_1^{M[G]},<)$;
\item  $(V^{N^*},<_V,P,E)$ is $(I^S,<,P,E)$;
\item  $W^{N^*}$ is $N_I$ (the $\tau$-functions and relations only
act on the $W$-sort);
\item  $R_1\subseteq U\times V$, with $R_1(\alpha,t)$ holding if and only if
$t\in J_\alpha$; and
\item $R_2\subseteq U\times W$, with $R_2(\alpha,b)$ holding if and only if
$b\in N_\alpha$.
\end{itemize}

%
%

Note that  $S\subseteq\omega_1^{M[G]}$ is a
$\tau^*$-definable subset of the $U$-sort of $N^*$ ($\alpha\in S$ if and only if $V\setminus R_1(\alpha,V)$ has a $<_V$-minimal element).  Also, on the
$W$-sort, the relation `$b\in\pcl(\abar)$' is definable by an
infinitary $\tau^*$-formula.  Thus, the relations `$b$ catches
$N_\alpha$' , `$b$ has unbounded reach in $N_\alpha$' and `$b$ has
bounded effect in $N_\alpha$' are each infinitarily $\tau^*$-definable
subsets of $U\times W$.

By construction, $N^*\models\psi$, where the infinitary $\psi$ asserts:
 `For every non-zero $\alpha\in U$, {\bf either} every element of $W^{N^*}$ that catches $N_\alpha$ also has unbounded reach in $N_\alpha$ {\bf or}
 there is an element of $W^{N^*}$ that catches $N_\alpha$ and has bounded effect in $N_\alpha$.'

To distinguish between these two possibilities,
 there is an infinitary $\tau^*$-formula $\theta(x)$ such that for $x$ from the $U$-sort, $\theta(x)$ holds if and only if there exists $b\in N_I\setminus N_{J_x}$
that catches and has unbounded reach in  $N_{J_x}$.
Thus, for non-zero $\alpha\in\omega_1^{M[G]}$ we have  $$N^*\models\theta(\alpha)\quad \Longleftrightarrow\quad \alpha\in S$$

Now, identify a countable fragment $L_\A$ of
$L_{\omega_1,\omega}(\tau^*)$ to include the formulas mentioned in
the last three paragraphs, along with infinitary formulas ensuring
$\tau$-atomicity.

Now, we switch our attention to $V$, and apply Theorem~\ref{transfer1} to
$(M[G],\epsilon)$, $L_\A$, and $N^*$.  This gives us a family $(M_X,E)$
of elementary extensions of $(M[G],\epsilon)$, each of size $\aleph_1$, indexed
by subsets $X\subseteq\omega_1$ ($=\omega_1^V$).  Each of these models of
ZFC has an $\tau^*$-structure, which we call $N^*_X$ inside it.  As well, for each $X\subseteq\omega_1$, there is a continuous,
strictly increasing mapping $t_X:\omega_1\rightarrow U^{N^*_X}$ with the property
that $$N^*_X\models\theta(t_X(\alpha))\quad\Longleftrightarrow\quad\alpha\in X$$
Let $(I^X,<^X,E^X,P^X)$ be the `$V$-sort' of $N^*_X$.  Clearly, each $I^X\in\II^*$.

Finally, the $W$-sort of each $\tau^*$-structure $N^*_X$ is the
universe of a $\tau$-structure, striated by $I^X$. We call this
`reduct' $N_X$. Note that by our choice of $L_\A$ and the fact that
$N^*_X\succeq_{L_\A} N^*$, we know that every $\tau$-structure $N_X$
is an atomic model of $T$ and is easily seen to be of cardinality
$\aleph_1$. Thus, the proof of Theorem~\ref{big} reduces to the
following:

\medskip\noindent{{\bf Claim.}}  If $X\setminus Y$ is stationary, then there is no
$\tau$-isomorphism $f:N_X\rightarrow N_Y$.
\medskip\par

\bp Fix $X,Y\subseteq\omega_1$ such that $X\setminus Y$ is stationary and by way of contradiction
assume that $f:N_X\rightarrow N_Y$ were a $\tau$-isomorphism.
%
%
Consider the
 $\tau^*$-structures $N^*_X$ and
$N^*_Y$ constructed above.  As notation, for each
$\alpha\in\omega_1^V$, let $N^X_\alpha$ and $N^Y_\alpha$ denote
$\tau$-elementary substructures with universes
$R_2(t_X(\alpha),N^*_X)$ and $R_2(t_Y(\alpha),N^*_Y)$, respectively.

Next, choose a club $C_0\subseteq\omega_1$ such that for every $\alpha\in C_0$:
\begin{itemize}
\item $\alpha$ is a limit ordinal;
\item  The restriction of $f:N_\alpha^X\rightarrow N_\alpha^Y$ is a $\tau$-isomorphism.
\end{itemize}

Denote the set of limit points of $C_0$ by $C$. As $C$ is club and
$(X\setminus Y)$ is stationary, choose $\alpha$ in their
intersection.
 Fix a strictly increasing $\omega$-sequence $\<\alpha_n:n\in\omega\>$ of elements from $C_0$ converging to $\alpha$.
As $\alpha\in X$, we can choose an element $b\in N_X\setminus N^X_\alpha$
 such that
$b$  catches $N^X_\alpha$ and has unbounded reach in $N^X_\alpha$.
That is, there is $\gamma<\alpha$ such that for every $\beta$
satisfying $\gamma<\beta<\alpha$,
$$\pcl(N^X_\gamma\cup\{b\},N_X)\cap N^X_\alpha\not\subseteq N^X_\beta.$$
Fix $n\in\omega$ such that $\alpha_n>\gamma$.  Then, for every $m\ge
n$
$$\pcl(N^X_{\alpha_n}\cup\{b\},N_Y)\cap N^X_\alpha\not\subseteq N^X_{\alpha_m}.$$

Thus, as `$b\in\pcl(\abar)$' is preserved under $\tau$-isomorphisms and $f[N^X_{\alpha_m}]=N^Y_{\alpha_m}$ setwise,
we have that $f(b)$ both catches and has unbounded reach in $N^Y_\alpha$.
As $\alpha\not\in Y$, we obtain a contradiction from  $N^*_Y\models\neg\theta(t_Y(\alpha))$ and $N^*_Y\models\psi$.

\bibliography{ssgroups}
\bibliographystyle{plain}

\end{document}